\documentclass[conference]{IEEEtran}
\IEEEoverridecommandlockouts
\usepackage{cite}
\usepackage{amsmath,amssymb,amsfonts}
\usepackage{algorithmic}
\usepackage{graphicx}
\usepackage{textcomp}
\usepackage{xcolor}
\usepackage[english]{babel}
\usepackage{amsthm}
\usepackage{url}
\usepackage{tikz}
\usepackage{mathabx}
\def\BibTeX{{\rm B\kern-.05em{\sc i\kern-.025em b}\kern-.08em
    T\kern-.1667em\lower.7ex\hbox{E}\kern-.125emX}}
\begin{document}

\title{QUBO formulation for the Snake-in-the-box and Coil-in-the-box problems\\
}

\author{\IEEEauthorblockN{Federico Fuidio}
\IEEEauthorblockA{\textit{Facultad de Ingeniería} \\
\textit{Universidad de Montevideo}\\
Montevideo, Uruguay \\
ffuidio@correo.um.edu.uy} 
\and
{Eduardo Canale}\\
\IEEEauthorblockA{\textit{Facultad de Ingeniería} \\
\textit{Universidad de la República}\\
Montevideo, Uruguay \\
canale@fing.edu.uy}
\and 
{Rafael Sotelo}\\
\IEEEauthorblockA{\textit{Facultad de Ingeniería} \\
\textit{Universidad de Montevideo}\\
Montevideo, Uruguay \\
rsotelo@um.edu.uy}
}

\maketitle

\begin{abstract}
This paper present the first QUBO formulations for the Snake-in-the-box (SITB) and Coil-in-the-box (CITB) problems. Both formulations are also capable of solving the NP-Hard problems of Maximum induced path and Maximum induced cylce respectively.  In the process we also found a new QUBO formulation for the Maximum Common Induced Sub-graph problem. We proved the correctness of our formulations for the SITB, CITB and Maximum Common Sub-graph problem, and tested the formulations of the SITB and CITB in both classical and quantum solvers, being able to get the best solution for up to 5 dimensions.
\end{abstract}


\section{Introduction}

\subsection{Snake-In-The-Box and Coil-In-The-Box problems}

Both the Snake-in-the-box problem (SITB problem for short) and the Coil-in-the-box problem (CITB problem for short) were first introduced by Kautz in 1958 \cite{Kautz_1958}. 

The SITB problem involves finding the longest induced path in the $n$-dimensional hypercube $Q_n$. (See section II for definitions of induced graphs and $Q_n$). In other words, it entails discovering the longest path that is isomorphic to an induced sub-graph of $Q_n$. Similarly, the CITB problem involves finding the longest induced cycle in the hypercube graph $Q_n$.

Note that in both cases, by \emph{longest} we mean the graph with the maximum number of edges. Since in a path the number of edges is one less than the number of vertices, and in a cycle both numbers are the same, both the SITB and CITB problems are equivalents if we aim to maximize the number of vertices instead than the number of edges. Therefore, in all the formulations presented here, we are trying to maximize the number of vertices. We call $n$-snake and $n$-coil to an induced path and cycle of $Q_{n}$ respectively. We often follow the notation introduced in \cite{Harary_1988} and call an $n$-snake and $n$-coil simply as a snake and a coil respectively.

Both problems have a wide range of applications, first, solutions of these problems can be used for error-detection codes \cite{Kautz_1958}. Also, they can be used for scalar quantizers in digital communication systems \cite{Kim_Neuhoff_2000}, rank modulation \cite{Yehezkeally_2012} or genetics \cite{Zinovik_2010}.

Usual methods to solve these problems involve brute force approaches \cite{Kautz_1958}, genetic algorithms \cite{Potter_1994} and Monte-Carlo Tree search \cite{Kinny_2012}. In table \ref{tab1}, we present the best known values for snake and coil lengths for dimensions 1 to 13. Notice that the length refer to the number of edges. Also, for dimentions 1 to 8, it is proved that the values shown are indeed the best ones.

\begin{table}[htbp]
\caption{Best knwon values for Snake and Coil in the box}
\begin{center}
\begin{tabular}{|c|c|c|}
\hline
\textbf{Dimension} & \textbf{\textit{SITB}}& \textbf{\textit{CITB}} \\
\hline
1 & \textbf{1}  & \textbf{0}  \\
2 & \textbf{2} \cite{Davies_1965} 	&	 \textbf{4} \cite{Kautz_1958} \\
3 & \textbf{4} \cite{Davies_1965}  	& 	 \textbf{6} \cite{Kautz_1958} \\
4 & \textbf{7} \cite{Davies_1965}  	& 	 \textbf{8} \cite{Kautz_1958} \\
5 & \textbf{13} \cite{Davies_1965} 	& 	 \textbf{14} \cite{Kautz_1958} \\
6 & \textbf{26} \cite{Davies_1965} 	& 	 \textbf{26} \cite{Davies_1965}  \\
7 & \textbf{50} \cite{Potter_1994} 	& 	 \textbf{48} \cite{Kochut_1996} \\
8 & \textbf{98} \cite{Ostergard_Pettersson_2015} 	& 	 \textbf{96} \cite{Ostergard_Pettersson_2014} \\
9 & 190 \cite{Wynn_2012} 			& 	188 \cite{Wynn_2012} \\
10 & 370 \cite{Kinny_2012}			&      366 \cite{Allison_Paulusma_2016} \\
11 & 712 \cite{Allison_Paulusma_2016} &       692 \cite{Allison_Paulusma_2016} \\
12 & 1373 \cite{Allison_Paulusma_2016} &     1344 \cite{Allison_Paulusma_2016} \\
13 & 2687 \cite{Allison_Paulusma_2016} &     2594 \cite{Allison_Paulusma_2016} \\
\hline
\multicolumn{3}{l}{\emph{Values in bold type are proven to be the best}}
\end{tabular}
\label{tab1}
\end{center}
\end{table}

\subsection{Quantum Annealer and QUBO formualtion}

QUBO (Quadratic Unconstrained Binary Optimization) \cite{Punnen_2022} is a mathematical model for optimization problems. It entails finding the minumum value of a binary quadratic form (i.e. finding the minimum value of a function of the form $\sum_{i, j} \alpha_{i, j}x_i x_j$ for $x_i \in \{0, 1\}$, $\alpha_{i, j} \in \mathbb{R}$). 

Quantum annealing \cite{Kadowaki_Nishimori_1998} represents a potential advantage in tackling the problem at hand compared to traditional algorithms. It is a method aimed at finding the global minimum of an objective function, drawing inspiration from the process of annealing in metallurgy where a material is heated and then slowly cooled to remove defects and attain a more stable structure.

D-Wave Systems \cite{D-Wave_Welcome} is a notable producer of commercial quantum annealers. These devices are specifically designed for quantum annealing and have found application in various research and commercial projects.

\subsection{Quantum Annealers and NP-Hard problems}

 NP-hard problems \cite{Garey1979} are distinguished by their substantial complexity, where the time required to find a solution increases exponentially with the size of the input. These problems are particularly challenging in computational theory due to their intractability using classical algorithms, especially as the problem size expands.

In classical computing, solving NP-hard problems typically involves heuristic or approximation methods, such as genetic algorithms \cite{Holland1992} or Monte-Carlo Tree search \cite{Browne2012}, to find near-optimal solutions within a reasonable timeframe. Exact solutions generally require brute-force approaches that are computationally expensive and impractical for large instances.

Quantum computing introduces a new paradigm for addressing NP-hard problems, leveraging the principles of quantum mechanics to process information in fundamentally different ways compared to classical computers. One promising approach is quantum annealing.

Quantum annealers offer potential advantages over classical solvers by exploring a vast solution space more efficiently due to quantum superposition and tunneling effects \cite{McGeoch2014}. Specifically, quantum annealers can escape local minima more effectively than classical methods, potentially finding better solutions for NP-hard problems.

\subsection{Structure of the paper}

The paper is organized as follows: Section II provides the necessary background in graph theory, introducing key concepts and definitions used throughout the paper. Section III introduces the QUBO formulations for the SITB and CITB problems, as well as formulations for the Induced Sub-graph problem and the Maximum Common Induced Sub-graph problem, which are foundational to the SITB and CITB formulations. Section IV presents the statements of the theorems related to the correctness of the QUBO formulations presented in section III. Section V discusses the results obtained from applying the QUBO formulations of the SITB and CITB problems using both classical and quantum solvers. Section VI concludes the paper, summarizing the findings and outlining future work, with a particular focus on exploring problem symmetries to reduce the number of variables in the QUBO formulations, aiming to solve larger instances of the problems. The appendix contains the detailed proofs of the theorems stated in Section IV, providing the mathematical rigor behind the correctness of the QUBO formulations.

\section{Background}

In this section we introduce some key concepts and definitions used in this paper.

\emph{Induced sub-graph: } Given a graph $G = (V, E)$, an induced subgraph by a subset of vertices $A \subseteq V$ is a subgraph of $G$ whose vertex set is equal to $A$, and whose edges consist of all edges connecting two vertices in $A$. We denote this graph as $G[A]$.

\emph{Graph Isomorphism: } We say that two graphs $G_1 = (V_1, E_1)$ and $G_2 = (V_2, E_2)$ are isomorphic (denoted as $G_1 \simeq G_2$) if there exists a bijective function $\phi : V_1 \rightarrow V_2$ that preserves the graph structure, i.e., $\{u, v \} \in E_1$ if and only if $\{\phi(u), \phi(v) \} \in E_2$. If such a function $\phi$ exists, we say that $\phi$ is an isomorphism between $G_1$ and $G_2$.

\emph{$n$-dimensional hypercube: }The hypecrube $Q_{n}$ can be defined as a graph whose vertex set is $V_{Q_{n}} = \{ 0, 1\}^{n}$ and edge set is $E_{Q_{n}} = \{ \{u, v \} \: | \: u, \: v \in V_{Q_{n}}, \: d_H(u, v) = 1 \}$, where $d_H(u, v)$ is the Hamming distance between $u$ and $v$.

\emph{Path and Coil graphs: } We define the graph $P_n$ as a path of $n$ vertices, meaning that $V_{P_n} = \{v_0, v_1, \dots, v_{n-1} \}$ and $E_{P_n} = \{ \{v_i, v_{i+1} \} \: | \: 0 \leq i < n-1\}$. Similarly, the graph $C_n$ is defined as a cycle of $n$ vertices, so that $V_{C_n} = \{v_0, v_1, \dots, v_{n-1} \}$ and $E_{C_n} = \{ \{v_i, v_{(i+1 \: \text{mod} \: n)} \} \: | \: 0 \leq i \leq n-1\}$.

\section{Formulation}

\subsection{QUBO formulation for the Induced sub-graph problem}

The Induced sub-graph problem \cite{Barrow_Burstall_1976} is a decision problem that involves determining if a graph $G_1 = (V_1, E_1)$ is an induced sub-graph of $G_2 = (V_2, E_2)$. Note that this problem is equivalent to determining if there exists an injective function $\phi : V_1 \rightarrow V_2$ such that $\{u, v\} \in E_1$ if and only if $\{\phi(u), \phi(v)\} \in E_2$. Calude et al. (2017) \cite{Calude_2017} provides a QUBO formulation for this problem. 

In this work, we present a similar formulation inspired by the formulation of the Graph Isomorphism Problem given by Lucas (2014) \cite{Lucas_2014}. The idea of the formulation involves defining binary variables to describe the function $\phi$. Let $x_{u, i}$ for $u \in V_1$ and $i \in V_2$ be such that $x_{u, i} = 1$ only if the vertex $u$ is mapped to $i$ by the function $\phi$. Define $s_i$ for $i \in V_2$ such that $s_i = 1$ only if there exists a vertex $u \in V_1$ that is mapped to $i$ by $\phi$. We arrive at the following QUBO formulation:

$$ Q = H_A + H_B \label{inducedGraph} $$

$$ H_A = \sum_{u \in V_1}\left( 1 - \sum_{i \in V_2} x_{u, i}\right)^2 + \sum_{i \in V_2}\left( s_i - \sum_{u \in V_1}x_{u, i}\right)^2\label{HAInducedGraph} $$

$$ H_B = \sum_{uv \in E_1}\sum_{ij \notin E_2} x_{u, i}x_{v, j} + \sum_{uv \notin E_1}\sum_{ij \in E_2} x_{u, i}x_{v, j}\label{HBInducedGraph} $$

This formulation needs a total of $|V_1||V_2| + |V_2|$ binary variables. 

We define $\phi$ as the relationship $\{(u, i) \in V_1 \times V_2 \: | \: x_{u, i} = 1 \}$.

Note that if $H_A = 0$, then $\phi$ defines an injective function from $V_1$ to $V_2$, and if $H_B = 0$, the function preserves the graph structure. Since $H_A \geq 0$ and $H_B \geq 0$, $\min(Q) \geq 0$, implying that if the minimum value of $Q$ is zero, then $H_A = H_B = 0$, and $G_1$ is an induced sub-graph of $G_2$.

\subsection{QUBO formulation for the Maximum Common Induced sub-graph problem}

The maximum common induced sub-graph problem \cite{Barrow_Burstall_1976} consists of finding, given two graphs $G_1$ and $G_2$, a graph with the maximum number of vertices that is an induced sub-graph of both $G_1$ and $G_2$. Equivalently, we aim to find the induced sub-graph of $G_1$ with the greatest order that is isomorphic to an induced sub-graph of $G_2$. \cite{Huang_2021} introduces the first QUBO formulation for this problem. In this paper, we provide a different formulation for this problem, that can be naturally extended to solve the SITB and CITB problems. 

We modify the formulation of the induced sub-graph problem by adding variables $p_u$ to indicate if a vertex $u$ in $V_1$ is mapped to a vertex in $V_2$. Let $x_{u, i}$ and $s_i$ be defined as in the previous formulation. Then, the QUBO formulation is as follows:

\begin{equation}
Q = \alpha H_A + \beta H_B + \gamma H_O\label{MaxInducedGraph}
\end{equation}

\begin{equation}
H_O = -\sum_{u \in V_1} p_{u}\label{HOMaxInducedGraph}
\end{equation}

\begin{equation}
H_A = \sum_{u \in V_1}\left( p_u - \sum_{i \in V_2} x_{u, i}	\right)^2 + \sum_{i \in V_2}\left( s_i - \sum_{u \in V_1} x_{u, i}\right)^2\label{HAMaxInducedGraph}
\end{equation}

\begin{equation}
H_B = \sum_{uv \in E_1}\sum_{ij \notin E_2} x_{u, i}x_{v, j} + \sum_{uv \notin E_1}\sum_{ij \in E_2} x_{u, i}x_{v, j}\label{HBMaxInducedGraph}
\end{equation}

Where $\alpha, \beta, \gamma > 0$. Also, we found that if $\alpha, \beta > \gamma$ then the minimum value of $Q$ is achieved only if we found a solution for the Maximum Common Induced Sub-graph problem.

This formulation uses $|V_1||V_2| + |V_1| + |V_2|$ binary variables.

As in the previous formulation, we define $\phi$ as the relationship $\{(u, i) \in V_1 \times V_2 \: | \: x_{u, i} = 1 \}$. Also, we define two sets $A_x \subseteq V_1$ and $B_x \subseteq V_2$ as:

\begin{equation}
A_{x} = \{ u \in V_1 \: | \: x_{u, i} = 1 \text{ for some } i \in V_2 \}
\end{equation}

\begin{equation}
B_{x} = \{ i \in V_2 \: | \: x_{u, i} = 1 \text{ for some } u \in V_1 \}
\end{equation}

In this formulation, $H_A = 0$ ensures that $\phi$ is an injective function from the set $A_x$ to $V_2$, and $H_B = 0$ ensures that the function preserves the graph structure. $H_O$ represents the objective function of the formulation, which aims to maximize the number of vertices in the induced sub-graph.

If we choose values for $\alpha$, $\beta$, $\gamma$ such that $\alpha, \beta > \gamma$, then the minimum value of $Q$ is achieved only if $G_1[A_x] \simeq G_2[B_x]$. We will prove this in the appendix.

\subsection{QUBO formulation for the SITB problem}

Using the formulation of the Maximum Common Induced sub-graph problem, we can solve both the SITB and CITB problems. The SITB problem can be viewed as an instance of the maximum common induced sub-graph problem, where $G_1 = P_{2^n}$ (the path of $2^n$ vertices) and $G_2 = Q_n$ (the $n$-dimensional hypercube graph). With an additional restriction to ensure that the graph $G_1[A_x]$ (the selected induced sub-graph of $P_n$) is connected.

In order to achieve this, we introduce a new term, $H_C$, to the previous formulation. The QUBO function $Q$ for the SITB problem is then defined as follows:

\begin{equation} \label{QUBO_SITB}
Q = \alpha H_A + \beta H_B + \gamma H_O + \delta H_C
\end{equation}

\begin{equation}
H_O = -\sum_{u \in V_1} p_{u}\label{HOSITB}
\end{equation}

\begin{equation}
H_A = \sum_{u \in V_1}\left( p_u - \sum_{i \in V_2} x_{u, i}	\right)^2 + \sum_{i \in V_2}\left( s_i - \sum_{u \in V_1} x_{u, i}\right)^2\label{HASITB}
\end{equation}

\begin{equation}
H_B = \sum_{uv \in E_1}\sum_{ij \notin E_2} x_{u, i}x_{v, j} + \sum_{uv \notin E_1}\sum_{ij \in E_2} x_{u, i}x_{v, j}\label{HBSITB}
\end{equation}

\begin{equation}
H_C = (1 - p_{v_0})^2 + \sum_{uv \in E_1} (p_u - p_v)^2 \label{HCSITB}
\end{equation}

Where $\alpha, \beta, \gamma, \delta > 0$, $\alpha, \beta > \gamma + 2 \delta$ and $\delta > 2^{n} \gamma$.

Since $|V_1| = |V_2| = 2^{n}$ this formulation uses $|V_1||V_2| + |V_1| + |V_2| =  2^{2n} + 2^{n+1}$ binary variables. 

Notice that if $H_C$ equals $1$, then the graph $G_1[A_x]$ is connected and contains the vertex $v_0$, and $H_C = 0$ if and only if all the variables $p_u$ are equal to one.

\subsection{QUBO formulation for the CITB problem}

The formulation for the Coil-in-the-box problem is not as straightforward. In the SITB problem, as we showed earlier, we can essentially solve the maximum common induced sub-graph problem by setting the largest path that can fit in the $n$-dimensional hypercube graph (a path with $2^n$ vertices) as $G_1$, setting $Q_n$ as $G_2$ and introducing an additional term to ensure that the chosen sub-graph of $P_{2^n}$ remains connected.

This approach solves the SITB problem for two reasons. First, any connected sub-graph of a path is also a path. Second, in a path, the term $H_C$ can ensure that the selected sub-graph of the path remains connected. However, it's important to note that in other types of graphs, such as general graphs or graphs with more complex structures, the $H_C$ term may not guarantee contentedness.

Then, in order to define the CITB problem, we cannot simply take $G_1$ as the largest possible cycle that can fit in the Hypercube graph (a cycle with $2^n$ vertices) and use the same formulation of the SITB. However we can slightly modify the SITB formulation in the following way:

First, define $G_1 = \left( V_1, \: E_1\right)$ where $V_1 = V_{\text{path}} \cup V_{\text{cycle}}$ and $E_1 = E_{\text{path}} \cup E_{\text{cycle}}$, and:

$$ V_{\text{path}} = \left\{ 1, 2, \dots, 2^n - 1 \right\} $$

$$ V_{\text{cycles}} = \left\{ (1,\: 3), \: (1, \: 4), \dots, \: (1, \: 2^n - 1)\right\} $$

$$ E_{\text{path}} = \left\{ \left\{j, \: j+1\right\}\: | \: 1 \leq j\leq 2^n - 2\right\} $$

$$ E_{\text{cycles}} =  \left\{ \left\{ i, v \right\}\: | \: i \in \{1, n\}, \:  v \in V_{\text{cycles}} \right\} $$

In Figure 1, we present a representation of the graph $G_1$ for $n = 3$. Vertices and edges belonging to $V_{\text{cycle}}$ and $E_{\text{cycle}}$ are drawn in dotted blue, while vertices and edges belonging to $V_{\text{path}}$ and $E_{\text{path}}$ are shown in red.

From $G_1$, we define an oriented graph $\widearrow{G_1}$, where the vertex set is $V_1$, and we assign a direction to each edge as shown in Figure 2. We refer to the edge set of the oriented graph as $\widearrow{E_1}$. Additionally, we denote the edge $(u, v) \in \widearrow{E_1}$ (with direction from $u$ to $v$) as $u \rightarrow v$.

\begin{figure} \label{Fig1}
\centering

\begin{tikzpicture}[main/.style = {draw, circle}]
\node[main, fill=red!20] (1) {$1$}; 
\node[main, fill=red!20] (2) [right of=1]{$2$}; 
\node[main, fill=red!20] (3) [right of=2] {$3$}; 
\node[main, fill=red!20] (4) [right of=3] {$4$}; 
\node[main, fill=red!20] (5) [right of=4] {$5$}; 
\node[main, fill=red!20] (6) [right of=5] {$6$}; 
\node[main, fill=red!20] (7) [right of=6] {$7$};
\node[main, fill=blue!20] (1_3) [above of=2] {$_{_{(1,3)}}$};
\node[main, fill=blue!20] (1_4) [above right of=1_3] {$_{_{(1,4)}}$};
\node[main, fill=blue!20] (1_5) [above right of=1_4] {$_{_{(1,5)}}$};
\node[main, fill=blue!20] (1_6) [above right of=1_5] {$_{_{(1,6)}}$};
\node[main, fill=blue!20] (1_7) [above right of=1_6] {$_{_{(1,7)}}$};
\draw[red, very thick] (1) -- (2); 
\draw[red, very thick] (2) -- (3);
\draw[red, very thick] (3) -- (4); 
\draw[red, very thick] (4) -- (5); 
\draw[red, very thick] (5) -- (6); 
\draw[red, very thick] (6) -- (7); 
\draw[dotted, blue, very thick] (1) to [out=90, in=180, looseness=1] (1_3);
\draw[dotted, blue, very thick] (1_3) to [out = 0, in = 90, looseness=1] (3);
\draw[dotted, blue, very thick] (1) to [out=90, in=180, looseness=1.2] (1_4); 
\draw[dotted, blue, very thick] (1_4) to [out=0, in=90, looseness=1] (4);
\draw[dotted, blue, very thick] (1) to [out=90, in=180, looseness=1.3] (1_5); 
\draw[dotted, blue, very thick] (1_5) to [out=0, in=90, looseness=1] (5);
\draw[dotted, blue, very thick] (1) to [out=90, in=180, looseness=1.4] (1_6); 
\draw[dotted, blue, very thick] (1_6) to [out=0, in=90, looseness=1] (6);
\draw[dotted, blue, very thick] (1) to [out=90, in=180, looseness=1.5] (1_7); 
\draw[dotted, blue, very thick] (1_7) to [out=0, in=90, looseness=1] (7);

\end{tikzpicture} \\
\caption{Representation of the graph $G_1$ for $n = 3$}

\end{figure}
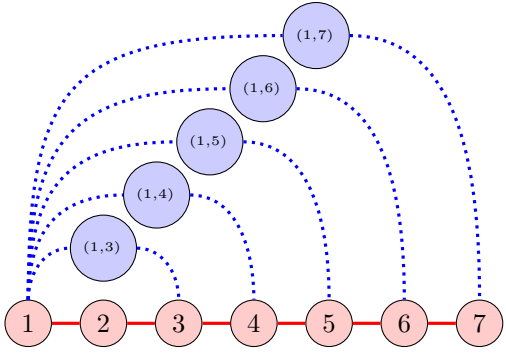

\begin{figure} \label{Fig2}
\centering

\begin{tikzpicture}[main/.style = {draw, circle}]
\node[main, fill=red!20] (1) {$1$}; 
\node[main, fill=red!20] (2) [right of=1]{$2$}; 
\node[main, fill=red!20] (3) [right of=2] {$3$}; 
\node[main, fill=red!20] (4) [right of=3] {$4$}; 
\node[main, fill=red!20] (5) [right of=4] {$5$}; 
\node[main, fill=red!20] (6) [right of=5] {$6$}; 
\node[main, fill=red!20] (7) [right of=6] {$7$}; 
\node[main, fill=blue!20] (1_3) [above of=2] {$_{_{(1,3)}}$};
\node[main, fill=blue!20] (1_4) [above right of=1_3] {$_{_{(1,4)}}$};
\node[main, fill=blue!20] (1_5) [above right of=1_4] {$_{_{(1,5)}}$};
\node[main, fill=blue!20] (1_6) [above right of=1_5] {$_{_{(1,6)}}$};
\node[main, fill=blue!20] (1_7) [above right of=1_6] {$_{_{(1,7)}}$};
\draw[red, very thick, ->] (1) -> (2); 
\draw[red, very thick, ->] (2) -> (3);
\draw[red, very thick, ->] (3) -> (4); 
\draw[red, very thick, ->] (4) -> (5); 
\draw[red, very thick, ->] (5) -> (6); 
\draw[red, very thick, ->] (6) -> (7); 
\draw[dotted, blue, very thick, <-] (1) to [out=90, in=180, looseness=1] (1_3);
\draw[dotted, blue, very thick, <-] (1_3) to [out = 0, in = 90, looseness=1] (3);
\draw[dotted, blue, very thick, <-] (1) to [out=90, in=180, looseness=1.2] (1_4); 
\draw[dotted, blue, very thick, <-] (1_4) to [out=0, in=90, looseness=1] (4);
\draw[dotted, blue, very thick] (1) to [out=90, in=180, looseness=1.3] (1_5); 
\draw[dotted, blue, very thick, <-] (1_5) to [out=0, in=90, looseness=1] (5);
\draw[dotted, blue, very thick, <-] (1) to [out=90, in=180, looseness=1.4] (1_6); 
\draw[dotted, blue, very thick, <-] (1_6) to [out=0, in=90, looseness=1] (6);
\draw[dotted, blue, very thick, <-] (1) to [out=90, in=180, looseness=1.5] (1_7); 
\draw[dotted, blue, very thick, <-] (1_7) to [out=0, in=90, looseness=1] (7);

\end{tikzpicture}\\
\caption{Representation of the directed graph $\widearrow{G_1}$ for $n = 3$}

\end{figure}
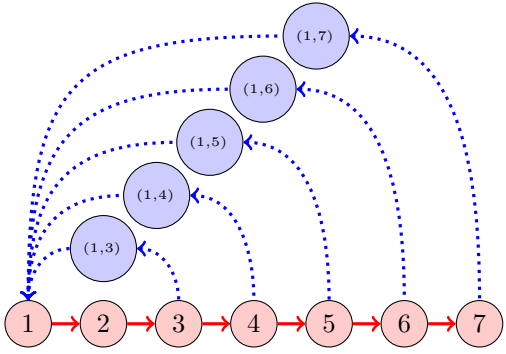

Finally, we define $G_2$ as the $n$-dimensional hypercube graph $Q_n$, and we define the binary variables $x_{u, i}$, $p_u$, and $s_i$ in the same manner as in our previous formulations. With these definitions, we arrive at the following QUBO formulation:

\begin{equation}
Q = \alpha H_A + \beta H_B + \gamma H_O + \delta H_C + \epsilon H_R\label{CITB}
\end{equation}

\begin{equation}
H_O = - \sum_{u \in E_1} p_u\label{HOCITB}
\end{equation}

\begin{equation}
H_A = \sum_{u \in V_1}\left( p_u - \sum_{i \in V_2} x_{u, i}	\right)^2 + \sum_{i \in V_2}\left( s_i - \sum_{u \in V_1} x_{u, i}\right)^2 \label{HACITB}
\end{equation}

\begin{equation}
H_{B} = \sum_{uv \in E_1}\sum_{ij \notin E_2} x_{u, i}x_{v, j} 	+ \sum_{uv \notin E_1}\sum_{ij \in E_2} x_{u, i}x_{v, j} \label{HBCITB}
\end{equation}

\begin{equation}
H_C =  \left(1 - \sum_{u \in V_{cycle}} p_u \right)^2 \label{HCCITB}
\end{equation}

\begin{equation}
H_R = \sum_{u \in V_{path}} \left( p_u - \sum_{v: \: u \rightarrow v \in \widearrow{E_1}}p_v \right)^2 \label{HRCITB}
\end{equation}

Where $\alpha, \beta, \gamma, \delta, \epsilon > 0$, $\alpha, \beta > \gamma + \delta + \epsilon$ and $\delta, \epsilon > 2^{n}\gamma$.

In this formulation $|V_1| = 2^{n+1} - 4$ and $|V_2| = 2^{n}$. Then the formulation uses $|V_1||V_2| + |V_1| + |V_2| = 2^{2n+1} - 2^{n} - 4$ binary variables.

Here, the term $H_C$, if null, ensures that exactly one vertex in $V_{\text{cycle}}$ is selected, representing the use of a single vertex from the cycle component. Whereas the term $H_R$, if null, ensures that for every selected vertex $u$ in $V_{\text{path}}$, there is exactly one corresponding vertex $v \in V_1$ such that $u \rightarrow v$ is in $\widearrow{E_1}$. Therefore, if both $H_C$ and $H_R$ are zero, the induced sub-graph $G_1[A_x]$ forms a cycle.

\subsection{Longest induced path and Longest induced cycle problems}

The longest induced path problem involves finding the longest induced path in an arbitrary graph $G$. Similarly, the longest induced cycle problem aims to find the longest induced cycle in an arbitrary graph $G$. Both of these problems are proven to be NP-hard \cite{Garey_Johnson_1990}.

Note that the SITB and CITB problems are just special cases of these problems when the graph $G$ is the hypercube $Q_n$. Thus, the formulations presented for the SITB and CITB problems can handle the more general NP-hard problems of the longest induced path and longest induced cycle in a graph, respectively.

Therefore, to solve the longest induced path (respectively, longest induced cycle) problem in a graph $G = (V, E)$, we can use the same formulation employed for the CITB (respectively, SITB) problem, where we set $G_2 = G$ instead of $G_2 = Q_n$. This demonstrates the versatility and broader applicability of our formulations in addressing these complex combinatorial problems.

\section{Correctness of the QUBO formulations}

In this section we state four theorems regarding the correctness of the formulations presented in this paper. The first two theorems (\ref{theorem_1} and \ref{theorem_2}) proves the correctness of the Maximum common induced sub-graph problem, under the assumption that $\alpha, \beta, \gamma > 0$ and $\alpha, \beta > \gamma$. These results are important since the formulation of both the SITB and CITB problems are based on our formulation for the Maximum common induced sub-graph problem.

Then we state two theorems (\ref{theorem_SITB} and \ref{theorem_CITB}) that proves the correctness of the QUBO formulations for the SITB and CITB problem respectively. The four theorems are proven in the Appendix.

\subsection{Correctness of the QUBO formulation of the Maximum Common Induced Sub-graph Problem}

\newtheorem{theorem}{Theorem}[section]

\begin{theorem} \label{theorem_1}
Let $Q(x, p, s)$ be the QUBO defined in \eqref{MaxInducedGraph} to \eqref{HBMaxInducedGraph} for the Maximum Common Induced Sub-graph problem, where $x, p, s$ are vectors containing the variables $x_{u, i}$, $p_u$ and $s_i$ respectively.

Then, if $\alpha, \: \beta > \gamma$ and $x, p, s$ are such that $Q(x, p, s)$ is the minimum value of Q then the relationship $\phi$ defined by $x$ is an injective function that preserves the graph structure.
\end{theorem}

\begin{theorem} \label{theorem_2}
(Correctness of the QUBO formulation for the Maximum Common Sub-graph Problem). Let $Q(x, p, s)$ be the QUBO defined in \eqref{MaxInducedGraph} to \eqref{HBMaxInducedGraph}, where $x, p, s$ are vectors containing the variables $x_{u, i}$, $p_u$ and $s_i$ respectively.

Then, if $\alpha, \: \beta > \gamma$ and $x, p, s$ are such that $Q(x, p, s)$ is the minimum value of Q then the relationship $\phi$ defined by $x$ is the solution of the Maximum Commom sub-graph problem.
\end{theorem}

\subsection{Correctness of the QUBO formulations for the SITB and CITB problems}

\begin{theorem} \label{theorem_SITB}
(Correctness of the QUBO formulation for the SITB Problem). Let $Q(x, p, s)$ be the QUBO defined in \eqref{QUBO_SITB} to \eqref{HCSITB} for the SITB problem, where $x, p, s$ are vectors containing the variables $x_{u, i}$, $p_u$ and $s_i$ respectively.

Then, if $\alpha, \: \beta > \gamma + 2\delta$, $\delta > 2^{n}\gamma$ and $x, p, s$ are such that $Q(x, p, s)$ is the minimum value of Q then the relationship $\phi$ defined by $x$ is an injective function that preserves the graph structure, and the subgraph $P_{2^{n}}[A_x]$ is a path.
\end{theorem}

\begin{theorem} \label{theorem_CITB}
(Correctness of the QUBO formulation for the CITB Problem). Let $Q(x, p, s)$ be the QUBO defined in \eqref{CITB} to \eqref{HRCITB} for the CITB problem, where $x, p, s$ are vectors containing the variables $x_{u, i}$, $p_u$ and $s_i$ respectively.

Then, if $\alpha, \: \beta > \gamma + \delta + \epsilon$, $\delta, \epsilon > 2^{n}\gamma$ and $x, p, s$ are such that $Q(x, p, s)$ is the minimum value of Q then the relationship $\phi$ defined by $x$ is an injective function that preserves the graph structure, and the subgraph $G_1[A_x]$ is a cycle.
\end{theorem}

\section{Results}

We tested both the SITB and CITB formulations in three different solvers. We use one classic solver and two QPU (Quantum Processing Unity) based solvers.

For the classic solver the SimulatedAnnealingSampler \cite{D-Wave_Simulated_Annealing} provided by D-Wave was used to simulate the implementation of a thermal annealing process. This simulator runs locally and was used on a personal laptop equipped with a 12th Gen Intel(R) Core(TM) i7-1255U CPU.

For the QPU based solvers we use two solvers that make use of the QPU. We used a Quantum-classical hybrid solver provided by D-Wave that uses both
classical and quantum resources \cite{D-Wave_Hybrid}.

For the purely QPU-based solution, we used the D-Wave Advantage\_system6.3 quantum computer, with 5614 qubits. We configure the number of shots to 1000, leaving the rest of the parameters unchanged. See \cite{D-Wave} (D-Wave documentation) for explanations for all parameters of the QPU. 

We accessed the QPU based solvers (both the hybrid and the purely quantum solver) using the \emph{sapi} interface \cite{Leap_Solvers}.

\subsection{Snake-in-the-box results}

Here we present the best solution obtained in the three solvers for the SITB problem. In both, the hybrid solver and the simulated annealer, we were able to get a $5$-snake of length 13 in $Q_5$. This is, we were able to get the best solution for dimension 5. For the QPU solver we only were able to get the best solution for dimension $3$.

Hybrid solution for $n = 5$ (note that it has $14$ vertices, therefore the length of the snake is $13$):
$$
\begin{matrix}
(10100) \rightarrow
(00100) \rightarrow
(00101) \rightarrow
(01101) \rightarrow\\
(01111) \rightarrow
(11111) \rightarrow
(10111) \rightarrow
(10011) \rightarrow\\
(00011) \rightarrow
(00010) \rightarrow
(01010) \rightarrow
(11010) \rightarrow\\
(11000) \rightarrow
(11001)
\end{matrix}
$$

Simulated solution for $n = 5$:
$$
\begin{matrix}
(11001) \rightarrow
(10001) \rightarrow
(10000) \rightarrow
(10100) \rightarrow\\
(00100) \rightarrow
(01100) \rightarrow
(01000) \rightarrow
(01010) \rightarrow\\
(01011) \rightarrow
(00011) \rightarrow
(00111) \rightarrow
(10111) \rightarrow\\
(11111) \rightarrow
(11110)
\end{matrix}
$$

QPU solution for $n = 3$:
$$
\begin{matrix}
(010)  \rightarrow
(000)  \rightarrow
(100)  \rightarrow
(101)  \rightarrow
(111)
\end{matrix}
$$

\subsection{Coil-in-the-box results}

For the CITB problem we were also able to get the best solution up to dimension 5 using the Simulated Annealing and the Hybrid Solver. For the QPU we only were able to get the best solution for dimension 2. 

Hybrid solution for $n = 5$:
$$
\begin{matrix}
(00101) \rightarrow
(00111) \rightarrow
(10111) \rightarrow
(11111) \rightarrow \\
(11101) \rightarrow
(11100) \rightarrow
(10100) \rightarrow
(10000) \rightarrow \\
(10010) \rightarrow
(11010) \rightarrow
(01010) \rightarrow
(01011) \rightarrow \\
(01001) \rightarrow
(00001)
\end{matrix}
$$

Simulated solution for $n = 5$:
$$
\begin{matrix}
(10111) \rightarrow
(10110) \rightarrow 
(00110) \rightarrow
(00100) \rightarrow \\
(01100) \rightarrow
(01000) \rightarrow
(11000) \rightarrow
(10000) \rightarrow \\
(10001) \rightarrow
(00001) \rightarrow
(00011) \rightarrow
(01011) \rightarrow \\
(11011) \rightarrow
(11111)
\end{matrix}
$$

QPU solution for $n = 2$:
$$
(11) \rightarrow
(01) \rightarrow
(00) \rightarrow
(10)
$$

\section{Conclusions and future work}

In this study, we introduced the first QUBO formulations for the Snake-in-the-box (SITB) and Coil-in-the-box (CITB) problems. Our formulations not only address these classic problems but also extend to solving the NP-hard problems of Maximum Induced Path and Maximum Induced Cycle, respectively. Through testing on both classical and quantum solvers, we demonstrated the efficacy of our approach by obtaining optimal solutions for dimensions up to 5.

The results obtained are promising, especially considering the complexity of the problems and the relatively early stage of quantum computing technology. The ability of our formulations to leverage quantum annealers for problem-solving hints at the potential for significant advancements in combinatorial optimization problems as quantum technology matures.

For future work, we aim to further optimize our QUBO formulations by exploiting the symmetries inherent in the SITB and CITB problems. By reducing the number of variables in the formulations, we hope to increase the scalability of our approach and enable the solving of these problems for higher dimensions on quantum annealers. This would not only advance the state-of-the-art for these specific problems but also contribute to the broader field of quantum optimization and its applications. Additionally, we plan to test our formulations using a wider variety of solvers and apply them to solve the more general problems of the longest induced path and longest induced cycle.

\newpage

\section{Appendix}

\subsection{Proof of theorem \ref{theorem_1}}

\textbf{\emph{Theorem \ref{theorem_1}}}: Let $Q(x, p, s)$ be the QUBO defined in \eqref{MaxInducedGraph} for the Maximum Induced Sub-graph problem, where $x, p, s$ are vectors containing the variables $x_{u, i}$, $p_u$ and $s_i$ respectively.

Then, if $\alpha, \: \beta > \gamma$ and $x, p, s$ are such that $Q(x, p, s)$ is the minimum value of Q then the relationship $\phi$ defined by $x$ is an injective function that preserves the graph structure.

\begin{proof}
We proceed by contradiction. First, suppose that the minimum value of $Q$ is achieved by $x$, $p$, $s$, and $\phi = \{(u, i) \in V_1 \times V_2 \: | \: x_{u, i} = 1 \}$ is either not a function from $A_x$ to $V_2$ (\emph{case 1}), not injective (\emph{case 2}), or does not preserve the graph structure (\emph{case 3}). Then, we find $\hat{x}$, $\hat{p}$, $\hat{s}$ with a strictly lower value of $Q$, leading to a contradiction.\\

\emph{Case 1: } The relationship $\phi$ is not a function from $A_x$ to $V_2$. This implies that there exists a vertex $u \in V_1$ that is mapped to more than one vertex in $V_2$, meaning that there exist $i, j \in V_2$ with $i \neq j$ such that $x_{u, i} = x_{u, j} = 1$. We can suppose that $p_u = s_i = s_j = 1$. This is because if any of these three variables were equal to 0, then we would immediately obtain a new solution with the corresponding variable equal to 1 and a strictly lower value of $Q$. Then, we have:

\begin{equation}
\sum_{k \in V_2} x_{u, k} \geq x_{u, i} + x_{u, j} = 2\label{case1_1}
\end{equation}

Using \eqref{case1_1} along with the fact that $p_u = 1$ we get:

\begin{equation}
p_u - \sum_{k \in V_2} x_{u, k} = 1 - \sum_{k \in V_2} x_{u, k} \leq -1 \label{case1_2}
\end{equation}

We can generate a new solution $\hat{x}, \hat{p}, \hat{s}$ with $\hat{x}_{u, j} = 0$, $\hat{s}_j = 0$, leaving all the other variables unchanged. We claim that this solution satisfy $Q(\hat{x}, \hat{p}, \hat{s}) < Q(x, p, s)$.

\begin{align*}
\sum_{k \in V_2} \hat{x}_{u, k} \: &= \: \hat{x}_{u, j} + \sum_{k \in V_2, \: k \not= j}\hat{x}_{u, k} \\
&= \: (x_{u, j} - 1) + \sum_{k\in V_2, \: k \not=j} x_{u, k} \\
&= \: \sum_{k \in V_2} x_{u, k} - 1
\end{align*}

Using now \eqref{case1_2}, along with the fact that $p_u = \hat{p_u}$:

$$\hat{p}_u - \sum_{k \in V_2}\hat{x}_{u, k} \:=\: p_u - \sum_{k \in V_2} x_{u, k} + 1 \: \leq \: 0$$
 
By taking absolute value on both sides we get:

$$\left| \hat{p}_{u} - \sum_{k \in V_2} \hat{x}_{u, k} \right| \:=\: \left| p_u - \sum_{k \in V_2} x_{u, k} \right| - 1$$

\begin{equation}
\left( \hat{p}_{u} - \sum_{k \in V_2} \hat{x}_{u, k} \right)^{2} \:<\: \left( p_u - \sum_{k \in V_2} x_{u, k} \right)^{2} \label{case1_3}
\end{equation}

The rest of the terms of $H_A$ (defined in \eqref{HAMaxInducedGraph}) remains unchanged. The term $H_O$ also remains unchanged, since $p_u = \hat{p}_u$. Finally, $H_B$ may change, but clearly it does not increase (since we are changing a variable that is set to 1 by 0). Therefore:

$$H_A(\hat{x}, \hat{p}, \hat{s}) \:<\: H_A(x, p, s)$$
$$H_B(\hat{x}, \hat{p}, \hat{s}) \: \leq \: H_B(x, p, s)$$
$$H_O(\hat{x}, \hat{p}, \hat{s}) \: =\: H_O(x, p, s)$$

Therefore $Q(\hat{x}, \hat{p}, \hat{s}) < Q(x, p, s)$, contradicting the minimality of $Q(x, p, s)$.\\

\emph{case 2: } The map $\phi$ is not injective, This implies that there exists a vertex $i \in V_2$ such that two (or more) vertices of $V_1$ are mapped to it. Meaning that there exists $i \in V_2$ and $u, v \in V_1$, $u \not= v$ such that $x_{u, i} = x_{v, i} = 1$. As in the first case, we can suppose that $s_i = p_u = p_v = 1$.

We take the solution $\hat{x}, \hat{p}, \hat{s}$ with $\hat{x}_{v, i} = 0$, $\hat{p}_{v} = 0$ and all the other values unchanged. We have then:

$$\beta H_B(\hat{x}, \hat{p}, \hat{s}) \: \leq \: \beta H_B(x, p, s)$$
$$\gamma H_O(\hat{x}, \hat{p}, \hat{s}) \:=\: \gamma H_O(x, p, s) + \gamma$$

With a similar analysis of the term $H_A$ that the one in case 1, we get:

\begin{equation}
\left( \hat{s}_i - \sum_{w \in V_2}\hat{x}_{w, i} \right)^2 \:<\: \left( s_i - \sum_{w \in V_2}x_{w, i} \right)^2
\end{equation}

Since we are working with natural numbers, we can  subtract one from the second term changing the strict inequality for a non-strict one.

 $$\left( \hat{s}_i - \sum_{w \in V_2}\hat{x}_{w, i} \right)^2 \: \leq \: \left( s_i - \sum_{w \in V_2}x_{w, i} \right)^2 - 1$$

All the other terms in $H_A$ remains unchanged. After multiplying by $\alpha$ we have:

$$\alpha H_A(\hat{x}, \hat{p}, \hat{s}) \: \leq \: \alpha H_A(x, p, s) - \alpha$$

Then:

$$Q(\hat{x}, \hat{p}, \hat{s}) \: \leq \: Q(x, p, s) + \gamma - \alpha$$

Since $\alpha > \gamma$, $Q(\hat{x}, \hat{p}, \hat{s}) < Q(x, p, s)$, contradicting the minimality of the solution.

Now, using what we proved in cases 1 and 2, we can assume that if $x, p, s$ gives the minimum value of $Q$, then $\phi$ is an injective function (implying that $H_A = 0$). This will be useful in the next case.\\

\emph{Case 3: } The function $\phi$ does not preserves the graph structure. Then, we can find $u, v \in V_1$, $i, j \in V_2$ such that $x_{u, i} = x_{v, j} = 1$, $\{u,v\} \in E_1$ but $\{i,j\} \notin E_2$ (The case when $\{u,v\} \notin E_1$ and $\{i,j\} \in E_2$ is analogous). Notice that this implies that $p_u = p_v = s_i = s_j = 1$.

Like we mentioned earlier, we can assume that $\phi$ is an injective function. In this case, we consider the solution $\hat{x}, \hat{p}, \hat{s}$ with $\hat{x}_{u, i} = 0$, $\hat{p}_u = 0$, and $\hat{s_i} = 0$. If $\phi$ is an injective function, then the map $\hat{\phi} = \{(u, i) \in V_1 \times V_2 \: | \: \hat{x}_{u, i}=1 \}$ is also an injective function, implying that $H_A$ remains equal to zero. Then:

$$\alpha H_A(\hat{x}, \hat{p}, \hat{s}) \:=\: \alpha H_A(x, p, s) \:=\: 0$$
$$\beta H_B(\hat{x}, \hat{p}, \hat{s}) \: \leq \: \beta H_B(x, p, s) - \beta $$
$$\gamma H_O(\hat{x}, \hat{p}, \hat{s}) \:=\: \gamma H_O(x, p, s) + \gamma$$

Implying that:

$$Q(\hat{x}, \hat{p}, \hat{s}) \: \leq \: Q(s, p, s) + \gamma - \beta$$

Since $\beta > \gamma$, $Q(\hat{x}, \hat{p}, \hat{s}) < Q(x, p, s)$, contradicting the minimality of $x, p, s$.\\

Thus, if the minimum value of $Q$ is achieved by the binary variables $x, p, s$ the set $\phi$ is an structure preserving injective function.

\end{proof}

\subsection{Proof of theorem \ref{theorem_2}}

\textbf{\emph{Theorem \ref{theorem_2}:}} (Correctness of the QUBO formulation for the Maximum Common Sub-graph Problem). Let $Q(x, p, s)$ be the QUBO defined in \eqref{MaxInducedGraph}, where $x, p, s$ are vectors containing the variables $x_{u, i}$, $p_u$ and $s_i$ respectively.

Then, if $\alpha, \: \beta > \gamma$ and $x, p, s$ are such that $Q(x, p, s)$ is the minimum value of Q then the relationship $\phi$ defined by $x$ is the solution of the Maximum Commom sub-graph problem.

\begin{proof}
Since $Q(x, p, s)$ is the minimum value of $Q$, using theorem \ref{theorem_1} we know that $\phi$ is an injective function from $A_x$ to $V_2$ that preserves the graph structure, meaning that $G_1[A_x] \simeq G_2[B_x]$. Then, we have $H_A(x, p, s) = H_B(x, p, s) = 0$ and $\gamma H_O(x, p, s) = - \gamma |A_x|$, therefore:

\begin{equation}
Q(x, p, s) = -\gamma |A_x| \label{Theorem2_eq1}
\end{equation}

By contradiction, suppose that there exists $A \subseteq V_1$, $B \subseteq V_2$ such that $G_1[A] \simeq G_2[B]$ with $|A| > |A_x|$.

Since $G_1[A] \simeq G_2[B]$, we know that there exists a bijective function $\psi : A \rightarrow B$ that preserves the graph structure. Define the binary variables $\hat{x}_{u, i}$, $\hat{p}_u$, $\hat{s}_i$ for $u \in V_1$, $i \in V_2$ in the following way:

$$
\hat{x}_{u, i}=
    \begin{cases}
        1 & \text{if }  u \in A \text{ and } \psi(u) = i \\
        0 & \text{otherwise} 
    \end{cases}
$$

$$
\hat{p}_u=
    \begin{cases}
        1 & \text{if }  u \in A \\
        0 & \text{otherwise} 
    \end{cases}
$$

$$
\hat{s}_i=
    \begin{cases}
        1 & \text{if }  i \in B \\
        0 & \text{otherwise} 
    \end{cases}
$$

The map $\{(u, i) \in V_1 \times V_2 \: | \: \hat{x}_{u, i} = 1\}$ is exactly $\psi$. This means that the map induced by the binary variables $\hat{x}_{u, i}$, $\hat{p}_u$ and $\hat{s}_{i}$ is an structure preserving injective function. Therefore $H_A(\hat{x}, \hat{p}, \hat{s}) = H_B(\hat{x}, \hat{p}, \hat{s}) = 0$ and $\gamma H_O(\hat{x}, \hat{p}, \hat{s}) = - \gamma |A|$, and:

$$Q(\hat{x}, \hat{p}, \hat{s}) = -\gamma |A|$$

This is a contradiction, since $Q(\hat{x}, \hat{p}, \hat{s}) < Q(x, p, s)$.

\end{proof}

\subsection{Proof of theorem of SITB}

\textbf{\emph{Theorem \ref{theorem_SITB}: }}(Correctness of the QUBO formulation for the SITB Problem). Let $Q(x, p, s)$ be the QUBO defined in \eqref{QUBO_SITB} to \eqref{HCSITB} for the SITB problem, where $x, p, s$ are vectors containing the variables $x_{u, i}$, $p_u$ and $s_i$ respectively.

Then, if $\alpha, \: \beta > \gamma + 2\delta$, $\delta > 2^{n}\gamma$ and $x, p, s$ are such that $Q(x, p, s)$ is the minimum value of Q then the relationship $\phi$ defined by $x$ is an injective function that preserves the graph structure, and the subgraph $P_{2^{n}}[A_x]$ is a path.

\begin{proof}
We can proceed exactly as in the proof of Theorem \ref{theorem_1}. First, we suppose that $Q(x, p, s)$ is the minimum value of $Q$ and $\phi$ is not a correct solution of the SITB problem, and then we find a new solution $\hat{x}, \hat{p}, \hat{s}$ with $Q(\hat{x}, \hat{p}, \hat{s}) < Q(x, p, s)$, leading to a contradiction. In this proof, we have 4 cases:\\

\emph{Case 1: } If the relationship $\phi$ is not a function, there exists a vertex $u \in V_1$ that is mapped to more than one vertex in $V_2$, meaning that there exist $i, j \in V_2$ with $i \neq j$ such that $x_{u, i} = x_{u, j} = 1$. Again, we can suppose that $p_u = s_i = s_j = 1$. With the same change of variables as the one

$$\alpha H_A(\hat{x}, \hat{p}, \hat{s}) < \alpha H_A(x, p, s)$$
$$\beta H_B(\hat{x}, \hat{p}, \hat{s}) \leq \beta H_B(x, p, s)$$
$$\gamma H_O(\hat{x}, \hat{p}, \hat{s}) = \gamma H_O(x, p, s)$$
$$\delta H_C(\hat{x}, \hat{p}, \hat{s}) = \delta H_C(x, p, s)$$

Implying that $Q(\hat{x}, \hat{p}, \hat{s}) < Q(x, p, s)$, contradicting the minimality of $Q(x, p, s)$.\\

\emph{Case 2: } $\phi$ is not injective. Again, we can do the same change as the one performed in the \emph{Case 2} of Theorem \ref{theorem_1}. Notice that in that after changing the value of only one $p_u$ then $H_C$ can change on at most $2$. Therefore we obtain:

$$\alpha H_A(\hat{x}, \hat{p}, \hat{s}) \leq \alpha H_A(x, p, s) - \alpha$$
$$\beta H_B(\hat{x}, \hat{p}, \hat{s}) \leq \beta H_B(x, p, s)$$
$$\gamma H_O(\hat{x}, \hat{p}, \hat{s}) = \gamma H_O(x, p, s) + \gamma$$
$$\delta H_C(\hat{x}, \hat{p}, \hat{s}) \leq \delta H_C(x, p, s) + 2\delta$$

Then:
 
$$Q(\hat{x}, \hat{p}, \hat{s}) \leq Q(x, p, s) + 2\delta + \gamma - \alpha$$

Since $\alpha > \gamma + 2\delta$ we have that $Q(\hat{x}, \hat{p}, \hat{s}) < Q(x, p, s)$, leading to a contradiction.\\

\emph{Case 3: } The function $\phi$ does not preserves the graph structure. Then we can perform the same change as in \emph{Case 3} of Theorem \ref{theorem_1}. In this case the effect on $H_C$ is the same, therefore we obtain:

$$H_A(\hat{x}, \hat{p}, \hat{s}) = H_A(x, p, s) = 0$$
$$\beta H_B(\hat{x}, \hat{p}, \hat{s}) \leq \beta H_B(x, p, s) - \beta$$
$$\gamma H_O(\hat{x}, \hat{p}, \hat{s}) = \gamma H_O(x, p, s) + \gamma$$
$$\delta H_C(\hat{x}, \hat{p}, \hat{s}) \leq \delta H_C(x, p, s) + 2\delta$$

Then:

$$Q(\hat{x}, \hat{p}, \hat{s}) \leq Q(x, p, s) + 2\delta + \gamma - \beta$$

Since $\beta > \gamma + 2\delta$ we have that $Q(\hat{x}, \hat{p}, \hat{s}) < Q(x, p, s)$, leading to a contradiction.\\

\emph{Case 4: } The selected subgraph of $P_{2^n}$ is not connected. Since $H_A(x, p, s) = H_B(x, p, s) = 0$, we have that $H_C(x, p, s) \geq 1$ (this is because $H_C(x, p, s) = 0$ only if all the $2^n$ vertices are used in the induced path, but we know that for $n \geq 2$ there is no induced path of $2^{n}$ vertices)

Since $G_1[A_x]$ is not connected then $H_C(x, p, s) > 1$. We can see that the trivial solution $\hat{x}_{u, i} = \hat{p}_u = \hat{s}_i = 0\ \forall u, i$  (all variables equal to zero) gives a lower value of $Q$:

$$Q(\hat{x}, \hat{p}, \hat{s}) = \delta$$
$$Q(x, p, s) \geq -|A_x|\gamma + 2\delta > -2^{n}\gamma + 2\delta$$

Therefore:

$$Q(x, p, s) \geq -2^{n}\gamma + Q(\hat{x}, \hat{p}, \hat{s}) + \delta$$

Since $\delta > 2^n\gamma$ we get that $Q(\hat{x}, \hat{p}, \hat{s}) < 0$, leading to a contradiction.

\end{proof}

\subsection{Proof of theorem CITB}

\textbf{\emph{Theorem \ref{theorem_CITB}: }} (Correctness of the QUBO formulation for the CITB Problem). Let $Q(x, p, s)$ be the QUBO defined in \eqref{CITB} to \eqref{HRCITB} for the CITB problem, where $x, p, s$ are vectors containing the variables $x_{u, i}$, $p_u$ and $s_i$ respectively.

Then, if $\alpha, \: \beta > \gamma + \delta + \epsilon$, $\delta, \epsilon > 2^{n}\gamma$ and $x, p, s$ are such that $Q(x, p, s)$ is the minimum value of Q then the map $\phi$ defined by $x$ is an injective function that preserves the graph structure, and the subgraph $G_1[A_x]$ cycle.

\begin{proof}
We proceed as in the proof of theorem $\ref{theorem_1}$. First, suppose that the minimum value of $Q$ is achieved by $x$, $p$, $s$, and $\phi = \{(u, i) \in V_1 \times V_2 \: | \: x_{u, i} = 1 \}$ is not a correct solution of the CITB problem, then, we find a solution $\hat{x}, \hat{p}, \hat{s}$ such that $Q(\hat{x}, \hat{p}, \hat{s}) < Q(x, p, s)$, leading to a contradiction. We have now 4 cases:\\

\emph{Case 1: } The relationship $\phi$ is not a function from $A_x$ to $V_2$. This implies that there exists a vertex $u \in V_1$ that is mapped to more than one vertex in $V_2$, meaning that there exist $i, j \in V_2$ with $i \neq j$ such that $x_{u, i} = x_{u, j} = 1$. Again, we can suppose that $p_u = s_i = s_j = 1$. Then we have:

$$\alpha H_A(\hat{x}, \hat{p}, \hat{s}) < \alpha H_A(x, p, s)$$
$$\beta H_B(\hat{x}, \hat{p}, \hat{s}) \leq \beta H_B(x, p, s)$$
$$\gamma H_O(\hat{x}, \hat{p}, \hat{s}) = \gamma H_O(x, p, s)$$
$$\delta H_C(\hat{x}, \hat{p}, \hat{s}) = \delta H_C(x, p, s)$$
$$\epsilon H_R(\hat{x}, \hat{p}, \hat{s}) \leq \epsilon H_R(x, p, s)$$

Then $Q(\hat{x}, \hat{p}, \hat{s}) < Q(x, p, s)$ contradicting the minimality of $Q(x, p, s)$.\\

\emph{case 2: } The map $\phi$ is not injective, This implies that there exists a vertex $i \in V_2$ such that two (or more) vertices of $V_1$ are mapped to it. Meaning that there exists $i \in V_2$ and $u, v \in V_1$, $u \not= v$ such that $x_{u, i} = x_{v, i} = 1$. As in the first case, we can suppose that $s_i = p_u = p_v = 1$.

We take the solution $\hat{x}, \hat{p}, \hat{s}$ with $\hat{x}_{v, i} = 0$, $\hat{p}_{v} = 0$ and all the other values unchanged. We have then:

$$\alpha H_A(\hat{x}, \hat{p}, \hat{s}) \leq \alpha H_A(x, p, s) - \alpha$$
$$\beta H_B(\hat{x}, \hat{p}, \hat{s}) \leq \beta H_B(x, p, s)$$
$$\gamma H_O(\hat{x}, \hat{p}, \hat{s}) = \gamma H_O(x, p, s) + \gamma$$
$$\delta H_C(\hat{x}, \hat{p}, \hat{s}) \leq \delta H_C(x, p, s) + \delta$$
$$\epsilon H_R(\hat{x}, \hat{p}, \hat{s}) \leq \epsilon H_R(x, p, s) + \epsilon$$

Then we have:

$$Q(\hat{x}, \hat{p}, \hat{s}) \leq Q(x, p, s) + \epsilon + \delta + \gamma - \alpha$$

Since $\alpha > \epsilon + \delta + \gamma$ we get $Q(\hat{x}, \hat{p}, \hat{s}) < Q(x, p, s)$, contradicting the minimality of $Q(x, p, s)$.\\

\emph{Case 3: } The function $\phi$ does not preserves the graph structure. Then, we can find $u, v \in V_1$, $i, j \in V_2$ such that $x_{u, i} = x_{v, j} = 1$, $\{u,v\} \in E_1$ but $\{i,j\} \notin E_2$ (The case when $\{u,v\} \notin E_1$ and $\{i,j\} \in E_2$ is analogous). Notice that this implies that $p_u = p_v = s_i = s_j = 1$.

Like we mentioned earlier, we can assume that $\phi$ is an injective function. In this case, we consider the solution $\hat{x}, \hat{p}, \hat{s}$ with $\hat{x}_{u, i} = 0$, $\hat{p}_u = 0$, and $\hat{s_i} = 0$. If $\phi$ is an injective function, then the map $\hat{\phi} = \{(u, i) \in V_1 \times V_2 \: | \: \hat{x}_{u, i}=1 \}$ is also an injective function, implying that $H_A$ remains equal to zero. Then:

$$\alpha H_A(\hat{x}, \hat{p}, \hat{s}) \:=\: \alpha H_A(x, p, s) \:=\: 0$$
$$\beta H_B(\hat{x}, \hat{p}, \hat{s}) \: \leq \: \beta H_B(x, p, s) - \beta $$
$$\gamma H_O(\hat{x}, \hat{p}, \hat{s}) \:=\: \gamma H_O(x, p, s) + \gamma$$
$$\delta H_C(\hat{x}, \hat{p}, \hat{s}) \leq \delta H_C(x, p, s) + \delta$$
$$\epsilon H_R(\hat{x}, \hat{p}, \hat{s}) \leq \epsilon H_R(x, p, s) + \epsilon$$

Then:

$$Q(\hat{x}, \hat{p}, \hat{s}) \leq Q(x, p, s) + \epsilon + \delta + \gamma - \beta$$

Since $\beta > \epsilon + \delta + \gamma$ we get $Q(\hat{x}, \hat{p}, \hat{s}) < Q(x, p, s)$, contradicting the minimality of $Q(x, p, s)$.\\

\emph{Case 4: } Finally, we prove that the selected induced sub-graph $G_1[A_x]$ is a cycle. For that we prove that $H_C(x, p, s) = H_R(x, p, s) = 0$. By contradiction, suppose that $H_C(x, p, s) > 0$, then, taking the trivial solution $\hat{x}_{u, i} = \hat{p}_{u} = \hat{s}_{i} = 0\ \forall u, i$.

$$Q(\hat{x}, \hat{p}, \hat{s}) = 0$$
$$Q(x, p, s) \geq -|A_x|\gamma + \delta + \epsilon H_R(x, p, s) > -2^{n}\gamma + \delta + \epsilon H_R(x, p, s)$$

Therefore, since $H_R(x, p, s) \geq 0$ and $\delta > 2^{n}\gamma$ we get that $Q(\hat{x}, \hat{p}, \hat{s}) < Q(x, p, s)$, leading to a contradiction.

Similarly, if $H_R(x, p, s) > 0$, we can take again the trivial solution $\hat{x}, \hat{p}, \hat{s}$:

$$Q(\hat{x}, \hat{p}, \hat{s}) = 0$$
$$Q(x, p, s) \geq -|A_x|\gamma + \epsilon + \delta H_C(x, p, s) > -2^{n}\gamma + \epsilon + \delta H_C(x, p, s)$$

Therefore, since $H_C(x, p, s) \geq 0$ and $\epsilon > 2^{n}\gamma$ we get that $Q(\hat{x}, \hat{p}, \hat{s}) < Q(x, p, s)$, leading to a contradiction.
\end{proof}

\end{document}